\documentclass[a4paper,notitlepage,11pt]{article}%
\usepackage{amssymb}
\usepackage{graphicx}
\usepackage{amsmath}
\usepackage{amsthm}
\usepackage{amsfonts}%
\providecommand{\U}[1]{\protect\rule{.1in}{.1in}}
\theoremstyle{plain}
\newtheorem{theorem}{Theorem}[section]

\newtheorem{corollary}[theorem]{Corollary}

\theoremstyle{definition}

\newtheorem{remark}[theorem]{Remark}
\numberwithin{equation}{section}
\numberwithin{theorem}{section}
\ifx\pdfoutput\relax\let\pdfoutput=\undefined\fi
\newcount\msipdfoutput
\ifx\pdfoutput\undefined\else
\ifcase\pdfoutput\else
\ifx\paperwidth\undefined\else
\ifdim\paperheight=0pt\relax\else\pdfpageheight\paperheight\fi
\ifdim\paperwidth=0pt\relax\else\pdfpagewidth\paperwidth\fi
\fi\fi\fi
\begin{document}

\title{Strictly positive solutions for one-dimensional nonlinear problems involving
the $p$-Laplacian\thanks{2000 \textit{Mathematics Subject Clasification}.
34B15; 34B18, 35J25, 35J61.} \thanks{\textit{Key words and phrases}. Elliptic
one-dimensional problems, indefinite nonlinearities, p-Laplacian, strictly
positive solutions.} \thanks{Partially supported by Secyt-UNC. The first
author would kindly like to dedicate this work to his teacher and friend
Tom\'{a}s Godoy.} }
\author{U. Kaufmann, I. Medri\thanks{\textit{E-mail addresses. }%
kaufmann@mate.uncor.edu (U. Kaufmann, Corresponding Author),
medri@mate.uncor.edu (I. Medri).}
\and \noindent\\{\small FaMAF, Universidad Nacional de C\'{o}rdoba, (5000) C\'{o}rdoba,
Argentina}}
\maketitle

\begin{abstract}
Let $\Omega$ be a bounded open interval, and let $p>1$ and $q\in\left(
0,p-1\right)  $. Let $m\in L^{p^{\prime}}\left(  \Omega\right)  $ and $0\leq
c\in L^{\infty}\left(  \Omega\right)  $. We study existence of strictly
positive solutions for elliptic problems of the form $-\left(  \left\vert
u^{\prime}\right\vert ^{p-2}u^{\prime}\right)  ^{\prime}+c\left(  x\right)
u^{p-1}=m\left(  x\right)  u^{q}$ in $\Omega$, $u=0$ on $\partial\Omega$. We
mention that our results are new even in the case $c\equiv0$.

\end{abstract}

\section{Introduction}

For $a<b$, let $\Omega:=(a,b)$, and let $p>1$ and $q\in\left(  0,p-1\right)
$. Let $m\in L^{p^{\prime}}\left(  \Omega\right)  $ and $0\leq c\in L^{\infty
}\left(  \Omega\right)  $. Our aim in this paper is to study the existence of
solutions for problems of the form%

\begin{equation}
\left\{
\begin{array}
[c]{ll}%
-\left(  \left\vert u^{\prime}\right\vert ^{p-2}u^{\prime}\right)  ^{\prime
}+c\left(  x\right)  u^{p-1}=m\left(  x\right)  u^{q} & \text{in }\Omega\\
u>0 & \text{in }\Omega\\
u=0 & \text{on }\partial\Omega.
\end{array}
\right.  \label{prob}%
\end{equation}
For applications we refer to \cite{branches} and the references therein.

When $c\equiv0$ and $0\not \equiv m\geq0$ it is known that (\ref{prob}) admits
a solution, see e.g. \cite{drabek}, Theorem 5.1, or \cite{chapa} and its
references for the case $p=2$. On the other hand, allowing $m$ to change sign
and under the assumption that $m\left(  x\right)  \geq m_{0}>0$ in some
$\overline{\Omega^{\prime}}\subset\Omega$, it can be proved that the problem
\begin{equation}
\left\{
\begin{array}
[c]{ll}%
-\left(  \left\vert u^{\prime}\right\vert ^{p-2}u^{\prime}\right)  ^{\prime
}=m\left(  x\right)  u^{q} & \text{in }\Omega\\
u=0 & \text{on }\partial\Omega
\end{array}
\right.  \label{casi}%
\end{equation}
possesses a \textit{nontrivial} \textit{nonnegative }solution (see Theorem 5.1
in \cite{drabek}, or \cite{orsina}, Section 5). We note however that in
general a (nontrivial) nonnegative solution of (\ref{casi}) need not be
strictly positive in $\Omega$ (in contrast to the superlinear case), and that
in fact the matter of existence of strictly positive solutions for these types
of problems is quite intriguing.

Recently, several \textit{non-comparable} sufficient conditions for the
existence of strictly positive solutions for (\ref{casi}) were exhibited in
\cite{ultimo} under some evenness assumptions on $m$ in the case $p=2$, and an
extension of some of these results for a (linear) strongly uniformly second
order elliptic operator was given in a paper \textquotedblleft Strictly
positive solutions for one-dimensional nonlinear elliptic
problems\textquotedblright, which has been submitted for publication by the
current authors. We refer to it later as [KM].

Let us mention that a natural way to attack these kind of problems is the well
known sub and supersolution method. Moreover, it is quite simple to provide
arbitrarily large supersolutions (see Remark \ref{homsup} below). In order to
construct the strictly positive subsolutions we shall adapt and extend the
approach developed in \cite{ultimo} and [KM]. Roughly speaking, we shall
divide $\Omega$ in parts, construct \textquotedblleft
subsolutions\textquotedblright\ en each of them and then find conditions on
$m$, $c$, $p$ and $q$ that guarantee that they can be joined accordingly to
obtain the desired subsolution. Certain conditions are presented in Theorem
\ref{bien}, and assuming that $m^{-}$ is essentially bounded further
non-comparable conditions are proved in Theorem \ref{aa} and Corollary
\ref{cerooo}.

Let us finally point out that although for the sake of simplicity we assume
that $c\geq0$, similar results can be obtained under some additional
assumptions if $c$ changes sign in $\Omega$ (see Remark \ref{negativo}).

\section{Preliminaries}

It is well known that for $g\in L^{1}\left(  \Omega\right)  $, the problem
$-\left(  \left\vert u^{\prime}\right\vert ^{p-2}u^{\prime}\right)  ^{\prime
}=g$ in $\Omega$, $u=0$ on $\partial\Omega$, admits a unique solution $u\in
C^{1}\left(  \overline{\Omega}\right)  $ such that $\left\vert u^{\prime
}\right\vert ^{p-2}u^{\prime}$ is absolutely continuous and that the equation
holds in the pointwise sense (e.g. \cite{ma1}, \cite{ma2}).

On the other side, it is also well known that if $g\in L^{p^{\prime}}(\Omega)$
(where as usual $p^{\prime}$ is given by $1/p+1/p^{\prime}=1$) and $0\leq c\in
L^{\infty}(\Omega)$, the problem%
\begin{equation}
\left\{
\begin{array}
[c]{ll}%
-\left(  \left\vert v^{\prime}\right\vert ^{p-2}v^{\prime}\right)  ^{\prime
}+c\left\vert v\right\vert ^{p-2}v=g & \text{in }\Omega\\
v=0 & \text{on }\partial\Omega
\end{array}
\right.  \label{g}%
\end{equation}
has a unique weak solution $v\in W_{0}^{1,p}\left(  \Omega\right)  $, i.e.,
satisfying
\[
\int_{\Omega}\left\vert v^{\prime}\right\vert ^{p-2}v^{\prime}\varphi^{\prime
}+c\left\vert v\right\vert ^{p-2}v\varphi=\int_{\Omega}g\varphi\qquad\text{for
all }\varphi\in W_{0}^{1,p}\left(  \Omega\right)
\]
(see e.g. \cite{garcia}). Furthermore, employing the comparison principles in
for instance \cite{les}, Chapter 6, and recalling the above paragraph, it is
easy to check that $v\in C^{1}\left(  \overline{\Omega}\right)  $, $\left\vert
v^{\prime}\right\vert ^{p-2}v^{\prime}$ is absolutely continuous and that
(\ref{g}) holds a.e. $x\in\Omega$.

We say that $0\leq v\in W_{0}^{1,p}\left(  \Omega\right)  $ is a (weak)
subsolution of (\ref{prob}) if
\begin{equation}
\int_{\Omega}\left\vert v^{\prime}\right\vert ^{p-2}v^{\prime}\varphi^{\prime
}+c\left(  x\right)  v^{p-1}\varphi\leq\int_{\Omega}m\left(  x\right)
u^{q}\varphi\qquad\text{for all }0\leq\varphi\in W_{0}^{1,p}\left(
\Omega\right)  \label{susu}%
\end{equation}
and $v=0$ on $\partial\Omega$; and $0\leq w\in W_{0}^{1,p}\left(
\Omega\right)  $ is said to be a supersolution if (\ref{susu}) holds (with $w$
in place of $v$) reversing the inequality, and $w\geq0$ on $\partial\Omega$.
The well known sub-supersolution method (\cite{drabek}, \cite{du}) gives a
solution provided there exist a subsolution $v$ and a supersolution $w$
satisfying $v\leq w$.

\begin{remark}
\label{homsup} Let us write as usual $m=m^{+}-m^{-}$ with $m^{+}=\max\left(
m,0\right)  $ and $m^{-}=\max\left(  -m,0\right)  $. If $m^{+}\not \equiv 0$,
one can readily verify that (\ref{prob}) admits arbitrarily large
supersolutions. Indeed, let $v\geq0$ be the solution of (\ref{g}) with $m^{+}$
in place of $g$, and let $k\geq(\left\Vert v\right\Vert _{\infty
}+1)^{q/(p-1-q)}$. Then $k(v+1)$ is a supersolution since $v=k>0$ on
$\partial\Omega$ and%
\begin{gather*}
-\left(  \left\vert k\left(  v+1\right)  ^{\prime}\right\vert ^{p-2}%
k(v+1)^{\prime}\right)  ^{\prime}+c\left(  k\left(  v+1\right)  \right)
^{p-1}\geq\\
k^{p-1}m^{+}\geq(k(\left\Vert v\right\Vert _{\infty}+1))^{q}m^{+}%
\geq(k(v+1))^{q}m\text{\qquad in }\Omega\text{. }\blacksquare
\end{gather*}

\end{remark}

The next remark summarizes some necessary facts about principal eigenvalues
for problems with weight involving the $p$-Laplacian operator.

\begin{remark}
\label{ppal}Let $0\leq c\in L^{\infty}(\Omega)$ and let $m\in L^{p^{\prime}%
}\left(  \Omega\right)  $ with $m^{+}\not \equiv 0$. There exists a positive
principal eigenvalue $\lambda_{1}\left(  m,\Omega\right)  $ and $\Phi\in
W_{0}^{1,p}\left(  \Omega\right)  $ satisfying
\begin{equation}
\left\{
\begin{array}
[c]{ll}%
-\left(  \left\vert \Phi^{\prime}\right\vert ^{p-2}\Phi^{\prime}\right)
^{\prime}+c\left(  x\right)  \Phi^{p-1}=\lambda_{1}\left(  m,\Omega\right)
m\left(  x\right)  \Phi^{p-1} & \text{in }\Omega\\
\Phi>0 & \text{in }\Omega\\
\Phi=0 & \text{on }\partial\Omega.
\end{array}
\right.  \label{auto}%
\end{equation}
Moreover, $\lambda_{1}\left(  m,\Omega\right)  $ is unique and simple (see
e.g. \cite{cuesta} and the references therein). $\blacksquare$
\end{remark}

\section{Main results}

In order to avoid overloading the notation, for $y\geq a$, $z\leq b$ and
$\varepsilon\geq0$ we set%
\[
M_{a,\varepsilon}^{-}\left(  y\right)  :=\int_{a}^{y}\left(  m^{-}\left(
x\right)  +\varepsilon\right)  dx,\qquad M_{b,\varepsilon}^{-}\left(
z\right)  :=\int_{z}^{b}\left(  m^{-}\left(  x\right)  +\varepsilon\right)
dx.
\]
If $\varepsilon=0$ we simply write $M_{a}^{-}\left(  y\right)  $ and
$M_{b}^{-}\left(  z\right)  $.

\begin{theorem}
\label{bien}Let $m\in L^{p^{\prime}}(\Omega)$ and suppose there exist $a\leq
x_{0}<x_{1}\leq b$ such that $0\not \equiv m\geq0$ in $I:=(x_{0},x_{1})$. Let%
\begin{gather}
\gamma:=\max\left\{  x_{1}-a,b-x_{0}\right\}  \qquad\text{and}\label{gamma}\\
\mathcal{M}_{p}:=\max{\large \{}M_{a}^{-}\left(  x_{1}\right)  ^{2-p}(%
{\textstyle\int\nolimits_{a}^{x_{1}}}
M_{a}^{-}\left(  x\right)  dx)^{p-1},M_{b}^{-}\left(  x_{0}\right)  ^{2-p}(%
{\textstyle\int\nolimits_{x_{0}}^{b}}
M_{b}^{-}\left(  x\right)  dx)^{p-1}{\large \}.}\nonumber
\end{gather}
(i) Assume $p\geq2$ and $q\in\left(  p-2,p-1\right)  $. If
\begin{gather}
\gamma^{p-2}\mathcal{M}_{2}<\frac{p-1}{\left(  p-1-q\right)  ^{p-1}}\frac
{1}{\lambda_{1}\left(  m,I\right)  }\qquad\text{and}\label{i1}\\
\gamma^{p}\left\Vert c\right\Vert _{L^{\infty}\left(  \Omega\right)  }%
\leq\frac{\left(  2-p+q\right)  \left(  p-1\right)  }{\left(  p-1-q\right)
^{p}}\label{i2}%
\end{gather}
then there \textit{exists a solution of }(\ref{prob})\textit{.}\newline(ii)
Assume $p\in\left(  1,2\right]  $. If%
\begin{gather}
\mathcal{M}_{p}<\frac{\left(  p-1\right)  ^{p}}{\left(  p-1-q\right)  ^{p-1}%
}\frac{1}{\lambda_{1}\left(  m,I\right)  }\qquad\text{and}\label{i3}\\
\gamma^{p}\left\Vert c\right\Vert _{L^{\infty}\left(  \Omega\right)  }%
\leq\left(  \frac{p-1}{p-1-q}\right)  ^{p}q\label{i4}%
\end{gather}
then there \textit{exists }a \textit{solution of }(\ref{prob})\textit{.}%
\newline
\end{theorem}

\textit{Proof}. Without loss of generality we assume that $a<x_{0}<x_{1}<b$
(in fact, it shall be clear from the proof how to proceed if either $x_{0}=a$
or $x_{1}=b$). Taking into account Remark \ref{homsup} it suffices to
construct a strictly positive (in $\Omega$) weak subsolution $u$ for
(\ref{prob}). Moreover, it is clear that it is enough to provide such
subsolution for (\ref{prob}) with $\tau m$ in place of $m$, for some $\tau>0$.

Let us prove (i). In view of (\ref{i1}) we may choose $\varepsilon>0$ small
enough and fix $\tau$ such that
\begin{equation}
\gamma^{p-2}\frac{\left(  p-1-q\right)  ^{p-1}}{p-1}\max\left\{  \int
_{a}^{x_{1}}M_{a,\varepsilon}^{-}\left(  x\right)  dx,\int_{x_{0}}%
^{b}M_{b,\varepsilon}^{-}\left(  x\right)  dx\right\}  \leq\frac{1}{\tau}%
\leq\frac{1}{\lambda_{1}\left(  m,I\right)  }.\label{tau1}%
\end{equation}
Let $x\in\left[  a,x_{1}\right]  $ and define%
\begin{gather}
u_{1}\left(  x\right)  :=\left(  \sigma\int_{a}^{x}M_{a,\varepsilon}%
^{-}\left(  y\right)  dy\right)  ^{k}\qquad\text{where}\label{def}\\
k:=\frac{1}{p-1-q},\qquad\sigma:=\frac{\tau\gamma^{p-2}}{\left(  p-1\right)
k^{p-1}}.\label{def2}%
\end{gather}
We have that $u_{1}(a)=0$ and that $u_{1}$ is strictly increasing. Also, from
the first inequality in (\ref{tau1}) it follows that $\left\Vert
u_{1}\right\Vert _{\infty}\leq1$. Let $l:=\left(  k-1\right)  \left(
p-1\right)  $. Since $q>p-2$ it holds that $l>0$. Furthermore,
\[
l-1+p=k\left(  p-1\right)  ,\qquad l+p-2=kq,
\]
and by (\ref{i2}) we also obtain that $k^{p-1}l\geq\gamma^{p}\left\Vert
c\right\Vert _{\infty}$. On the other hand, since $M_{a,\varepsilon}^{-}$ is
strictly increasing we derive that $\left(  x-a\right)  M_{a,\varepsilon}%
^{-}\left(  x\right)  \geq\int_{a}^{x}M_{a,\varepsilon}^{-}\left(  y\right)
dy$ for all $x$. Taking into account the aforementioned facts, (\ref{def2})
and that $p\geq2$ and $x_{1}-a\leq\gamma$, some computations show that
\begin{gather}
-\left(  \left\vert u_{1}^{\prime}\left(  x\right)  \right\vert ^{p-2}%
u_{1}^{\prime}\left(  x\right)  \right)  ^{\prime}=\label{cuenta}\\
-\left(  k\sigma^{k}\right)  ^{p-1}\left(  l\left(  \int_{a}^{x}%
M_{a,\varepsilon}^{-}\left(  y\right)  dy\right)  ^{l-1}M_{a,\varepsilon}%
^{-}\left(  x\right)  ^{p}+\right.  \nonumber\\
\left.  \left(  p-1\right)  \left(  \int_{a}^{x}M_{a,\varepsilon}^{-}\left(
y\right)  dy\right)  ^{l}M_{a,\varepsilon}^{-}\left(  x\right)  ^{p-2}\left(
m^{-}\left(  x\right)  +\varepsilon\right)  \right)  \leq\nonumber\\
-\left(  k\sigma^{k}\right)  ^{p-1}\left(  \frac{l}{\gamma^{p}}\left(
\int_{a}^{x}M_{a,\varepsilon}^{-}\left(  y\right)  dy\right)  ^{k\left(
p-1\right)  }+\right.  \nonumber\\
\left.  \frac{\left(  p-1\right)  }{\gamma^{p-2}}\left(  \int_{a}%
^{x}M_{a,\varepsilon}^{-}\left(  y\right)  dy\right)  ^{kq}m^{-}\right)
\leq\nonumber\\
-\left\Vert c\right\Vert _{\infty}\sigma^{k\left(  p-1\right)  }\left(
\int_{a}^{x}M_{a,\varepsilon}^{-}\left(  y\right)  dy\right)  ^{k\left(
p-1\right)  }-\tau m^{-}\sigma^{kq}\left(  \int_{a}^{x}M_{a,\varepsilon}%
^{-}\left(  y\right)  dy\right)  ^{kq}\leq\nonumber\\
-cu_{1}^{p-1}-\tau m^{-}u_{1}^{q}\leq-cu_{1}^{p-1}+\tau mu_{1}^{q}\text{\qquad
in }\left(  a,x_{1}\right)  .\nonumber
\end{gather}

In a similar way, if for $x\in\left[  x_{0},b\right]  $ we set $u_{3}$ by
$u_{3}\left(  x\right)  :=\left(  \sigma\int_{x}^{b}M_{b,\varepsilon}%
^{-}\left(  y\right)  dy\right)  ^{k}$ with $k$ and $\sigma$ given by
(\ref{def2}), then $u_{3}(b)=0$, $u_{3}$ is strictly decreasing, $\left\Vert
u_{3}\right\Vert _{\infty}\leq1$ and
\[
-\left(  \left\vert u_{3}^{\prime}\right\vert ^{p-2}u_{3}^{\prime}\right)
^{\prime}+cu_{3}^{p-1}\leq\tau mu_{3}^{q}\text{\qquad in }\left(
x_{0},b\right)  .
\]

On the other side, let $u_{2}>0$ with $\left\Vert u_{2}\right\Vert
_{L^{\infty}\left(  I\right)  }=1$ be the positive principal eigenfunction
associated to the weight $m$ in $I$, that is satisfying (\ref{auto}) with $I$
in place of $\Omega$. Recalling that $m\geq0$ in $I$ and that $q<p-1$, from
the second inequality in (\ref{tau1}) we get%
\[
-\left(  \left\vert u_{2}^{\prime}\right\vert ^{p-2}u_{2}^{\prime}\right)
^{\prime}+cu_{2}^{p-1}=\lambda_{1}\left(  m,I\right)  mu_{2}^{p-1}\leq\tau
mu_{2}^{q}\text{\qquad in }I.
\]

Since $u_{1}\left(  a\right)  =u_{3}\left(  b\right)  =u_{2}\left(
x_{0}\right)  =u_{2}\left(  x_{1}\right)  =0$ and $\left\Vert u_{1}\right\Vert
_{\infty},\left\Vert u_{3}\right\Vert _{\infty}\leq1=\left\Vert u_{2}%
\right\Vert _{\infty}$, arguing as in the proof of Theorem 3.1 (i) in [KM] we
can find $\underline{x}_{0},\overline{x}_{1}\in\Omega$ with $\underline{x}%
_{0}<\overline{x}_{1}$ and such that
\begin{gather}
u_{1}(\underline{x}_{0})=u_{2}(\underline{x}_{0}),\text{\qquad}u_{2}\left(
\overline{x}_{1}\right)  =u_{3}\left(  \overline{x}_{1}\right)  ,\label{ollo}%
\\
u_{1}^{\prime}(\underline{x}_{0})\leq u_{2}^{\prime}(\underline{x}%
_{0}),\text{\qquad}u_{2}^{\prime}(\overline{x}_{1})\leq u_{3}^{\prime
}(\overline{x}_{1}).\nonumber
\end{gather}
We now define a function $u$ by $u:=u_{1}$ in $\left[  a,\underline{x}%
_{0}\right]  $, $u:=u_{2}$ in $\left[  \underline{x}_{0},\overline{x}%
_{1}\right]  $ and $u:=u_{3}$ in $\left[  \overline{x}_{1},b\right]  $. Taking
into account (\ref{ollo}), a simple integration by parts yields that $u$ is a
weak subsolution for (\ref{prob}) with $\tau m$ in place of $m$, and as we
said at the beginning of the proof this proves (i).

Let us prove (ii). We first pick $\varepsilon>0$ sufficiently small and take
$\tau$ such that%
\begin{align}
\frac{\left(  p-1-q\right)  ^{p-1}}{\left(  p-1\right)  ^{p}}\left(  \int
_{a}^{x_{1}}M_{a,\varepsilon}^{-}\left(  x\right)  dx\right)  ^{p-1} &
\leq\frac{1}{\tau}\leq\frac{1}{M_{a,\varepsilon}^{-}\left(  x_{1}\right)
^{2-p}\lambda_{1}\left(  m\right)  },\label{tau2}\\
\frac{\left(  p-1-q\right)  ^{p-1}}{\left(  p-1\right)  ^{p}}\left(
\int_{x_{0}}^{b}M_{b,\varepsilon}^{-}\left(  x\right)  dx\right)  ^{p-1} &
\leq\frac{1}{\tau}\leq\frac{1}{M_{b,\varepsilon}^{-}\left(  x_{0}\right)
^{2-p}\lambda_{1}\left(  m\right)  }.\label{tau3}%
\end{align}
(this is possible due to (\ref{i3})) Let $M_{\varepsilon}:=\max\left\{
M_{a,\varepsilon}^{-}\left(  x_{1}\right)  ,M_{b,\varepsilon}^{-}\left(
x_{0}\right)  \right\}  $. We shall build a strictly positive subsolution for
(\ref{prob}) with $\tau M_{\varepsilon}^{p-2}m$ in place of $m$. For
$x\in\left[  a,x_{1}\right]  $ we set $u_{1}$ as in (\ref{def}) with
\[
k:=\frac{p-1}{p-1-q},\qquad\sigma:=\frac{1}{k}\left(  \frac{\tau}{p-1}\right)
^{1/(p-1)}%
\]
in place of (\ref{def2}). Again $u_{1}(a)=0$, $u_{1}$ is strictly increasing
and using the first inequality in (\ref{tau2}) one can check that $\left\Vert
u_{1}\right\Vert _{\infty}\leq1$. Taking $l$ as in (i) we now obtain $l=kq$
and also as before we have $l-1+p=k\left(  p-1\right)  $ and $k^{p-1}%
l\geq\gamma^{p}\left\Vert c\right\Vert _{\infty}$. Furthermore, recalling that
$p\leq2$ and arguing as in (\ref{cuenta}) we deduce that
\begin{gather*}
-\left(  \left\vert u_{1}^{\prime}\left(  x\right)  \right\vert ^{p-2}%
u_{1}^{\prime}\left(  x\right)  \right)  ^{\prime}\leq-\left(  k\sigma
^{k}\right)  ^{p-1}\left(  l\left(  \int_{a}^{x}M_{a,\varepsilon}^{-}\left(
y\right)  dy\right)  ^{l-1+p}+\right.  \\
\left.  \left(  p-1\right)  \left(  \int_{a}^{x}M_{a,\varepsilon}^{-}\left(
y\right)  dy\right)  ^{l}M_{a,\varepsilon}^{-}\left(  x_{1}\right)
^{p-2}m^{-}\right.  \leq\\
-\left(  k\sigma^{k}\right)  ^{p-1}\left(  \frac{l}{\gamma^{p}}\left(
\int_{a}^{x}M_{a,\varepsilon}^{-}\left(  y\right)  dy\right)  ^{k\left(
p-1\right)  }+\right.  \\
\left.  \left(  p-1\right)  \left(  \int_{a}^{x}M_{a,\varepsilon}^{-}\left(
y\right)  dy\right)  ^{kq}M_{\varepsilon}^{p-2}m^{-}\right)  \leq\\
-\left\Vert c\right\Vert _{\infty}\sigma^{k\left(  p-1\right)  }\left(
\int_{a}^{x}M_{a,\varepsilon}^{-}\left(  y\right)  dy\right)  ^{k\left(
p-1\right)  }-\tau M_{\varepsilon}^{p-2}m^{-}\sigma^{kq}\left(  \int_{a}%
^{x}M_{a,\varepsilon}^{-}\left(  y\right)  dy\right)  ^{kq}\leq\\
-cu_{1}^{p-1}-\tau M_{\varepsilon}^{p-2}m^{-}u_{1}^{q}\leq-cu_{1}^{p-1}+\tau
M_{\varepsilon}^{p-2}mu_{1}^{q}\text{\qquad in }\left(  a,x_{1}\right)
\text{.}%
\end{gather*}
Since $u_{3}$ can be defined analogously and, taking into account the
definition of $M_{\varepsilon}$ and the second inequality in (\ref{tau2}) and
(\ref{tau3}), $u_{2}$ can be chosen as above (i.e. as the normalized positive
principal eigenfunction with respect to the weight $m$ in $I$), reasoning as
in (i) the theorem follows. $\blacksquare$

\begin{remark}
\label{integral}(i) Let us note that when $m\in C\left(  \Omega\right)  $ the
condition $0\not \equiv m\geq0$ in $I$ is necessary in order to have a
(nontrivial) nonnegative solution for (\ref{prob}). \newline(ii) Let us also
observe that if $p=2$ then (\ref{i1})-(\ref{i2}) coincide with (\ref{i3}%
)-(\ref{i4}) and that the resulting conditions extend the ones in [KM] (see
Theorem 3.5 (ii) there). $\blacksquare$
\end{remark}

For $p>1$ and $q\in\left(  0,p-1\right)  $ we set
\begin{equation}
C_{p,q}:=\left(  \frac{p}{p-1-q}\right)  ^{p-1}\frac{\left(  p-1\right)
\left(  q+1\right)  }{p-1-q}. \label{cpq}%
\end{equation}
We point out that for any $p>1$, $\lim_{q\rightarrow p-1^{-}}C_{p,q}=\infty$.
We shall now assume that $m^{-}\in L^{\infty}(\Omega)$. In the following
theorem we suppose that $c\not \equiv 0$, the case $c\equiv0$ is considered in
Corollary \ref{cerooo} below.

\begin{theorem}
\label{aa}Assume $c\not \equiv 0$. Let $m\in L^{p^{\prime}}(\Omega)$ with
$m^{-}\in L^{\infty}(\Omega)$ and suppose there exist $a\leq x_{0}<x_{1}\leq
b$ such that $0\not \equiv m\geq0$ in $I:=(x_{0},x_{1})$. Let $\gamma$ and
$C_{p,q}$ be given by (\ref{gamma}) and (\ref{cpq}) respectively.\newline(i)
Assume $p\geq2$. If
\begin{equation}
\frac{\left\Vert m^{-}\right\Vert _{L^{\infty}(\Omega)}}{\left\Vert
c\right\Vert _{L^{\infty}(\Omega)}}\sinh^{p}\left(  \left(  \frac{\left\Vert
c\right\Vert _{L^{\infty}(\Omega)}}{C_{p,q}}\right)  ^{1/p}\gamma\right)
\leq\frac{1}{\lambda_{1}(m,I)}\label{seno}%
\end{equation}
then there \textit{exists a} \textit{solution of }(\ref{prob})\textit{.}%
\newline(ii) Assume $p\in\left(  1,2\right)  $. If%
\begin{equation}
\frac{\left\Vert m^{-}\right\Vert _{L^{\infty}(\Omega)}}{\left\Vert
c\right\Vert _{L^{\infty}(\Omega)}}\left(  e^{\left(  \frac{\left\Vert
c\right\Vert _{L^{\infty}(\Omega)}}{C_{p,q}}\right)  ^{1/p}\gamma}-1\right)
^{p}\leq\frac{1}{\lambda_{1}(m,I)}\label{expo}%
\end{equation}
then there \textit{exists a} \textit{solution of }(\ref{prob})\textit{.}%
\newline
\end{theorem}

\textit{Proof. }The proof follows the lines of the proof of Theorem \ref{aa}
and hence we omit the details. We only indicate briefly how to construct
$u_{1}$ in both (i) and (ii). Suppose first (\ref{seno}) holds. Let $\tau$ be
such that
\begin{equation}
\frac{\left\Vert m^{-}\right\Vert _{L^{\infty}(\Omega)}}{\left\Vert
c\right\Vert _{L^{\infty}(\Omega)}}\sinh^{p}\left(  \left(  \frac{\left\Vert
c\right\Vert _{L^{\infty}(\Omega)}}{C_{p,q}}\right)  ^{1/p}\gamma\right)
\leq\frac{1}{\tau}\leq\frac{1}{\lambda_{1}(m,I)}\label{taun}%
\end{equation}
and for $x\in\left[  a,x_{1}\right]  $ define
\[
f(x):=\left(  \frac{\tau\left\Vert m^{-}\right\Vert _{\infty}}{\left\Vert
c\right\Vert _{\infty}}\right)  ^{1/p}\sinh\left(  \left(  \frac{\left\Vert
c\right\Vert _{\infty}}{C_{p,q}}\right)  ^{1/p}\left(  x-a\right)  \right)  .
\]
It is easy to check that $\left(  C_{p,q}^{1/p}f^{\prime}\right)  ^{2}-\left(
\left\Vert c\right\Vert _{\infty}^{1/p}f\right)  ^{2}=\left(  \tau\left\Vert
m^{-}\right\Vert _{\infty}\right)  ^{2/p}$ in $(a,x_{1})$. Moreover, $f(a)=0$,
$f$ is increasing (in particular, employing (\ref{taun}) and the fact that
$x_{1}-a\leq\gamma$ we see that $\left\Vert f\right\Vert _{\infty}\leq1$) and
$f^{\prime\prime}\geq0$ in $\left(  a,x_{1}\right)  $. Let us now choose
\begin{equation}
k:=\frac{p}{p-1-q},\qquad l:=\left(  k-1\right)  \left(  p-1\right)
.\label{kl}%
\end{equation}
It holds that $l-1=kq$, $l-1+p=k\left(  p-1\right)  $ and $k^{p-1}l=C_{p,q}$.
Define $u_{1}:=f^{k}$. Taking into account the above mentioned facts and that
$p\geq2$ we find that in $\left(  a,x_{1}\right)  $
\begin{gather}
-\left(  \left\vert u_{1}^{\prime}\right\vert ^{p-2}u_{1}^{\prime}\right)
^{\prime}+cu^{p-1}\leq\label{teo2}\\
-k^{p-1}\left(  lf^{l-1}\left(  f^{\prime}\right)  ^{p}+\left(  p-1\right)
f^{l}\left(  f^{\prime}\right)  ^{p-2}f^{\prime\prime}\right)  +\left\Vert
c\right\Vert _{\infty}f^{k\left(  p-1\right)  }\leq\nonumber\\
-k^{p-1}lf^{l-1}\left(  f^{\prime}\right)  ^{p}+\left\Vert c\right\Vert
_{\infty}f^{k\left(  p-1\right)  }=-f^{l-1}\left(  C_{p,q}\left(  f^{\prime
}\right)  ^{p}-\left\Vert c\right\Vert _{\infty}f^{p}\right)  \leq\nonumber\\
-f^{l-1}\left(  \left(  C_{p,q}^{1/p}f^{\prime}\right)  ^{2}-\left(
\left\Vert c\right\Vert _{\infty}^{1/p}f\right)  ^{2}\right)  ^{p/2}%
=-f^{l-1}\tau\left\Vert m^{-}\right\Vert _{\infty}\leq\tau mu_{1}^{q}%
\text{.}\nonumber
\end{gather}

Suppose now (\ref{expo}) holds. In this case we take $\tau$ and $f$ such that%
\begin{gather*}
\frac{\left\Vert m^{-}\right\Vert _{L^{\infty}(\Omega)}}{\left\Vert
c\right\Vert _{L^{\infty}(\Omega)}}\left(  e^{\left(  \frac{\left\Vert
c\right\Vert _{L^{\infty}(\Omega)}}{C_{p,q}}\right)  ^{1/p}\gamma}-1\right)
^{p}\leq\frac{1}{\tau}\leq\frac{1}{\lambda_{1}(m,I)},\\
f\left(  x\right)  :=\sigma\left(  e^{\lambda\left(  x-a\right)  }-1\right)
\text{,\qquad where}\\
\sigma:=\left(  \frac{\tau\left\Vert m^{-}\right\Vert _{\infty}}{\left\Vert
c\right\Vert _{\infty}}\right)  ^{1/p}\text{,\qquad}\lambda:=\left(
\frac{\left\Vert c\right\Vert _{\infty}}{C_{p,q}}\right)  ^{1/p}.
\end{gather*}
Let $k$ and $l$ be given by (\ref{kl}), and let $u_{1}:=f^{k}$. Reasoning as
in (\ref{teo2}) yields
\begin{gather*}
-\left(  \left\vert u_{1}^{\prime}\right\vert ^{p-2}u_{1}^{\prime}\right)
^{\prime}+cu^{p-1}\leq-f^{l-1}\left(  C_{p,q}\left(  f^{\prime}\right)
^{p}-\left\Vert c\right\Vert _{\infty}f^{p}\right)  =\\
-f^{l-1}\left(  C_{p,q}\left(  \sigma\lambda\right)  ^{p}e^{p\lambda\left(
x-a\right)  }-\left\Vert c\right\Vert _{\infty}\sigma^{p}\left(
e^{\lambda\left(  x-a\right)  }-1\right)  ^{p}\right)  \leq\\
-f^{l-1}\left\Vert c\right\Vert _{\infty}\sigma^{p}=-f^{l-1}\tau\left\Vert
m^{-}\right\Vert _{\infty}\leq\tau mu_{1}^{q}\text{\qquad in }\left(
a,x_{1}\right)  .\text{ }\blacksquare
\end{gather*}

\qquad

\begin{remark}
\label{pri}A quick look of the proof of the above theorem shows that (ii)
holds for any $p>1$. We observe however that one can verify that the
inequality (\ref{seno}) is better than (\ref{expo}). $\blacksquare$
\end{remark}

\begin{corollary}
\label{cerooo} Let $m$ be as in the above theorem and suppose $c\equiv0$. If
\begin{equation}
\frac{\left\Vert m^{-}\right\Vert _{L^{\infty}(\Omega)}\gamma^{p}}{C_{p,q}%
}\leq\frac{1}{\lambda_{1}(m,I)} \label{cortito}%
\end{equation}
then there \textit{exists a} \textit{solution of }(\ref{prob})\textit{.}%
\newline
\end{corollary}

\textit{Proof. }It is enough to note that the left side of either (\ref{seno})
or (\ref{expo}) tend to the left side of (\ref{cortito}) when $\left\Vert
c\right\Vert _{\infty}$ goes to zero. $\blacksquare$

\begin{remark}
\label{negativo} Let us suppose that $c$ changes sign in $\Omega$. An
inspection of the proofs of the theorems shows that one can still argue in the
same way as before replacing $c$ by $c^{+}$ in order to construct the
functions $u_{1}$ and $u_{3}$. Furthermore, if the positive principal
eigenvalue $\lambda_{1}\left(  m,I\right)  $ exists (for necessary and
sufficient conditions on this question, see \cite{cuesta}, Section 2) and if
the problem (\ref{g}) with $m^{+}$ in place of $g$ admits a nonnegative
solution, then all the analogous results to the case $c\geq0$ can be proven
allowing $c$ to change sign in $\Omega$. $\blacksquare$
\end{remark}

\end{document}